# Estimating Relevant Portion of Stability Region using Lyapunov Approach and Sum of Squares

Chetan Mishra, James S. Thorp, Virgilio A. Centeno
Bradley Dept. of Electrical and Computer Engineering,
Virginia Polytechnic Institute and State University,
Blacksburg-24061, Virginia, USA
Email: {chetan31, jsthorp, virgilio}@vt.edu

Anamitra Pal
School of Electrical, Computer, and Energy Engineering,
Arizona State University, Tempe-85281, Arizona, USA
Email: Anamitra.Pal@asu.edu

*Abstract*—Traditional Lyapunov-based transient stability assessment approaches focus on identifying the stability region (SR) of the equilibrium point under study. When trying to estimate this region using Lyapunov functions, the shape of the final estimate is often limited by the degree of the function chosen – a limitation that results in conservativeness in the estimate of the SR. More conservative the estimate is in a particular region in state space, smaller is the estimate of the critical clearing time (CCT) for disturbances that drive the system towards that region. In order to reduce this conservativeness, we propose a methodology that uses the disturbance trajectory data to skew the shape of the final Lyapunov-based SR estimate. We exploit the advances made in the theory of sum of squares decomposition to algorithmically estimate this region. The effectiveness of this technique is demonstrated on a power systems classical model.

## I. Introduction

The transient stability analysis of power systems traditionally relied on time domain simulations that were computationally intensive, thereby making online applications challenging. This led to the ongoing development of various direct methods [1] which start by defining a scalar function whose specific level set approximates the stability region (SR) for the post fault system configuration. The direct methods only required knowledge of the system state values at the time of fault clearing to predict the stability, thus significantly reducing the simulation time. Direct methods were of two types, namely, analytical energy function-based and Lyapunov-based. The techniques belonging to the first category start by characterizing the stability boundary and then defining an analytical energy function that is proved to exist under certain system related assumptions. This function is constrained to decrease along all the post-fault system trajectories and can take both positive and negative values. Using a unique level set of the energy function for estimating the SR (also known as the closest unstable equilibrium point approach) was seen to be extremely conservative in terms of the size of the final SR estimate relative to the actual SR size [1]. This led to the idea of relevant stability boundary and the controlling unstable equilibrium point (CUEP) approach where only the portion of the SR relevant to the particular disturbance was estimated using a particular level set of the energy function thus reducing the conservativeness. This chosen level set differed according to the disturbance being studied. While these approaches relied on a derived analytical energy function proven to satisfy the required criteria, there were numerous system related assumptions (such as transverse intersection of stable and unstable manifolds on the stability boundary, hyperbolicity and finiteness of equilibrium points, low transmission conductances, etc. [1]) that were impossible to verify and/or were found to be violated under various system conditions [2].

The Lyapunov-based techniques aim at estimating a scalar function referred to as the Lyapunov function [3] whose maximal level set represents an estimate of the SR. One advantage of these techniques is that they do not rely on any assumptions made regarding the underlying system dynamics. However, unlike the previous category of techniques, the Lyapunov function as well as its maximal level set needs to be estimated each time for changing system and is not necessarily parametric. Traditionally, this was achieved based on experience and/or domain specific knowledge and was a challenging task. Zubov [4] proposed a technique for algorithmically estimating a Lyapunov function that gave an exact SR estimate. The effectiveness of this technique was shown for power systems in [5]. However, the technique revolved around solving a partial differential equation which did not, in general, possess a closed form solution and therefore became impractical to solve for even medium-sized systems. With the advances made in the area of positive polynomials, sum of squares (SOS) programming [6] has become a powerful tool for analyzing complex polynomial dynamical systems [7]. An algorithmic construction of the Lyapunov function and the associated SR estimate using SOS programming referred to as the expanding interior algorithm was proposed in [8] and was applied to power systems classical model in [9]. This approach was further extended to do stability analysis of power systems structure preserving models by Anghel et al. [10]. Furthermore, the decomposition approaches proposed in [11] were also utilized to analyze larger scale interconnected systems. The increased deployment of renewables in the grid has added a new twist to this already challenging problem. The tripping of these generators resulting from ride through characteristics transformed the power system from a purely autonomous non-linear dynamical system to a constrained system [12]. Lyapunov based approaches have been fairly efficient in dealing with such systems while the analytical energy function based methods haven't caught up further justifying the need for these approaches.

Now, the Lyapunov's direct method aims at finding the largest estimate of the SR which is followed by assessing the stability for a given disturbance trajectory. However, large is a vague term making it difficult to compare to different estimates. A more effective way of comparing two different SR estimates is the amount of time the disturbance trajectory is contained inside each which is basically how well an estimate estimates the relevant portion of the SR. Therefore, in this paper, we try incorporating the information from the disturbance trajectory in the SR estimation process to reduce the conservativeness in the CCT estimate.

The rest of the paper is structured as follows. The estimation of SR using SOS for a classical power system model is presented in Section II. The idea of relevant portion of SR is discussed in Section III with a methodology proposed to choose the shape of the expanding region for SR estimate in the algorithm in [9]. Effectiveness of the proposed methodology is demonstrated in Section IV through time-domain simulations. The conclusion is provided in Section V.

## II. SUM OF SQUARES AND APPLICATION TO TSA

### A. Stability Region Estimation for Power Systems Classical Model using Expanding Interior Algorithm

For a system $\dot{x} = f(x)$ with an equilibrium at the origin, given a positive scalar function $V(x)$ with $V(0) = 0$ and a region $D \subseteq \mathbb{R}^n$ containing the equilibrium such that $\dot{V}(x) < 0 \; \forall x \neq 0 \in D$, any region $\Omega_\beta := \{x \in D | V(x) \leq \beta\}$ is a positive invariant region contained in the equilibrium's SR [3]. Thus, for a given Lyapunov function $V$, its largest level set contained inside the region $D$ gives an estimate of the SR. Therefore, the aim is to find a function whose corresponding maximal level set covers the biggest portion of the actual SR. However, power systems is not a polynomial system and SOS programming is valid only for polynomial systems. The power systems classical model with uniform damping in single machine reference frame is shown below.

$$\dot{\delta}_{in_g} = \omega_{in_g} \quad (1)$$

$$M_i \dot{\omega}_{in_g} = P_{m_i} - \left(\sum_{j=1:n_g} E_i E_j \left(G_{ij} \cos\left(\delta_{in_g} - \delta_{jn_g}\right) + B_{ij} \sin\left(\delta_{in_g} - \delta_{jn_g}\right)\right)\right) + -\frac{M_i}{M_{n_g}}(P_{m_{n_g}} - \sum_{j=1:n_g} E_{n_g} E_j \left(G_{n_g j} \cos\left(-\delta_{jn_g}\right) + B_{n_g j} \sin\left(-\delta_{jn_g}\right)\right)\right) - D_i \omega_{in_g} \quad (2)$$

$$\delta_{in_g} = \delta_i - \delta_{n_g}, \omega_{in_g} = \omega_i - \omega_{n_g} \quad (3)$$
$$i = 1,2 \ldots (n_g - 1)$$

The variable transformation shown in Table I [9] is used to convert the above system into a DAE polynomial system of the form: $\dot{z} = f(z)$ and $0 = g(z)$, where $z$ denotes the states of the polynomial system.

**Table II: Variable transformation for classical model**

| New Variable | Function of Original States |
|---|---|
| $z_i$ | $\omega_{in_g} - \omega_{in_g}^s(0) = \omega_{in_g}$ |
| $z_{n_g+2i-2}$ | $\sin(\delta_{in_g} - \delta_{in_g}^s)$ |
| $z_{n_g+2i-1}$ | $1 - \cos(\delta_{in_g} - \delta_{in_g}^s)$ |

The SR estimate for the system created using Table III can be obtained by utilizing the expanding interior algorithm [8]. The idea is to expand a semi-algebraic set $P_\beta = \{z | p(z) \leq \beta\}$ that is forced to be contained inside the level set of an unknown Lyapunov function $V(z)$ that satisfies the conditions discussed in [3]. The algorithm starts with an initial local estimate of the SR which can be obtained by solving the following SOS optimization problem.

$$\max_{V, \dot{V}, s_2, s_6, \lambda_1, \lambda_2} \beta_1 \quad (4)$$
$$-s_2(\beta_1 - p_1) + V - \lambda_1^T g - l_1 \text{ is SOS}$$
$$-s_6(\beta_1 - p_1) - \dot{V} - \lambda_2^T g - l_2 \text{ is SOS}$$

The shape of $p_1$ determines the function space the initial Lyapunov is to be found in as well as the invariant set defined by it. It may not be possible to find a decently sized invariant set contained inside a randomly shaped $p_1 \leq \beta_1$. A safe choice for $p_1$ that is usually used is a sphere described by $z^T z$.

$$\max_{s_1, s_2, s_3, \lambda_1} c \quad (5)$$
$$-s_1(-c+V) - s_2(p_1 - \beta_1) + s_3(-c+V)(p_1 - \beta_1)$$
$$- \lambda_1 g - (p_1 - \beta_1)^2 \text{ is SOS}$$
$$V = \frac{V}{c}, \dot{V} = \frac{\dot{V}}{c}$$

The expanding interior algorithm can be written as [8]:

$$\max \beta \quad (6)$$
$$s_2 V - \lambda_1^T g - l_1 \text{ is SOS}$$
$$-s_6(\beta - p) - \lambda_2^T g - (V-1) \text{ is SOS}$$
$$-s_8(1-V) - s_9 \dot{V} - \lambda_3^T g \text{ is SOS}$$

It was mentioned in [9] that a proper choice of $p(z)$ was necessary to get a good estimate and a methodology was proposed to replace $p(z)$ by the Lyapunov function $V$ obtained using that $p(z)$ and re-running the algorithm. This step was repeated till $P_\beta$ converged to $\{V \leq c\}$. However, it was noticed by us that the starting value of $p(z)$ heavily influenced the shape of the final estimate. The reason being that each successive Lyapunov estimate was forced to contain the previous ones. This can be seen in Fig. 1 for the SR estimates obtained for Example A of [9] for different choices of initial value of $p(z)$. The blue estimate is from a spherical $p(z)$ which is usually the default choice. The red estimate is obtained from an elliptical $p(z)$ that has its major axis along the direction of machine 2 states. This reflects in the shape of the final estimate which is longer in the $\delta_2$ direction. Similarly the green estimate is obtained when the major axis is along the direction of machine 1 states. Here it should be kept in mind that we are trying to estimate the SR in terms of the polynomial states while the final estimate in Fig. 1 is plotted in the original system states.

Past research has not addressed the criterion for choosing the initial value of the function $p(z)$. Usually the default choice for this function is $z^T z$ (sphere). As mentioned earlier, since our aim is TSA for one or more disturbance trajectories, a careful starting choice of $p(z)$ could help expand the Lyapunov estimate closer to the part of the actual stability boundary those trajectories are heading towards. This leads us to solving the challenging problem of (i) constructing $p(z)$ from the disturbance trajectory data and (ii) finding the relevant portion of the SR.

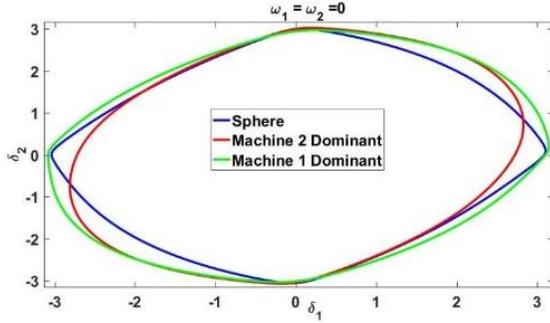

Fig. 1: SR estimate for different starting value of $p(z)$

### III. RELEVANT PORTION OF STABILITY REGION BASED METHODOLOGY

#### A. Relevant Portion of Stability Region (SR)

During a typical disturbance, the power system goes through pre-disturbance, disturbance-on and post-disturbance phases with transient stability normally estimated for the post-disturbance phase. Consider the two machine system whose actual SR is shown in Fig. 2. Two different Lyapunov estimates in $\omega_1 = \omega_2 = 0$ plane along with the fault trajectory is also shown. From the figure, it can be observed that the fault has a localized effect on machine 1, resulting in a more horizontal fault trajectory. Though the Lyapunov estimate $V_1$ gives a better estimate of the overall SR as compared to $V_2$, it is not so for analyzing this particular disturbance trajectory. It can be seen that longer the disturbance sustains, closer the system is pushed towards the stability boundary. The question then becomes: *how long one can sustain the disturbance without going outside the SR of the post-disturbance system?* A good SR estimate in this regard would be one that contains the disturbance trajectory for the longest time. Therefore, for this particular fault trajectory, the estimate from $V_2$ will be better than that from $V_1$, and hence the critical clearing time (CCT) estimate obtained using $V_2$ will be more accurate than that obtained using $V_1$.

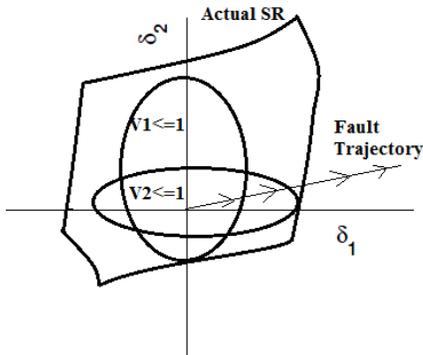

Fig. 2: Two machine system with different Lyapunov SR estimate

As mentioned before, this is the governing idea behind the CUEP approach where the unstable equilibrium point on the portion of the stability boundary that the disturbance trajectory is trying to head towards is found with the constant energy surface passing through it serving as the relevant SR. There is however a fundamental difference in what can be achieved when trying to estimate the relevant portion of SR through a single Lyapunov function vs. the CUEP approach. While we may be able to choose along what axes the Lyapunov estimate must be expanded, it is tricky to also incorporate the direction in the expansion process. For example, in Fig. 2, the disturbance trajectory heads towards positive and not negative $\delta_1$ direction but the default choice of $p(z)$ is not able to take that into consideration. This can be overcome by choosing complex shapes for $p(z)$. However, it would make the strategy even more intensive, from a computational perspective. Thus, there is a genuine need to calculate a suitable shaped $p(z)$ algorithmically for a given disturbance trajectory data.

#### B. Principal Component Analysis (PCA) Based Approach to Choosing $P_\beta$ using Fault Trajectory Data

As mentioned before, looking for higher degree Lyapunov functions for multi-machine systems is not practical due to its computational complexity. Also, it has been shown in [9] that quadratic Lyapunov functions give a good enough SR estimate. Thus, in this paper, we try to find a positive quadratic function $p(z)$ belonging to a family of ellipsoids. A good candidate for an ellipsoid would be one whose axes are aligned in the directions of the maximum variances of the disturbance trajectory dataset with the axis lengths reflecting the variance in the trajectory data in each direction. The successive level sets of such an ellipsoid would expand less in the directions where the disturbance trajectories do not head. A simple candidate for this ellipsoid can be obtained by applying principal component analysis (PCA) [13] to the raw disturbance data as opposed to the traditional approach of centering and normalizing the data. The reason for not using centered data is that the center of the final ellipse has to be at the origin, which is the location of the stable equilibrium point (SEP). Also, normalizing it will diminish the relative axes lengths and result in an ellipse that expands with regards to the absolute magnitude changes in every state variable. We know that in larger power systems, the effect of the disturbances that threaten the transient stability (mainly faults) are localized. Therefore, it is difficult to create disturbances that result in fault trajectories which cause significant displacement along all the states. This implies that almost always, the ellipse needed would have largely differing axes lengths.

Essentially, PCA finds a set of orthogonal vectors aligned in the direction of variances of the provided data. Also provided are the eigenvalues associated with each orthogonal vector which provide information about the relative variances in the data captured by each direction. The proposed candidate for initial $p(z)$ is as follows,

$$p(z) = z^T A z \tag{7}$$

$$A = Emat \times \begin{pmatrix} \frac{1}{\sqrt{\lambda_1}} & \cdots & 0 \\ \vdots & \ddots & \vdots \\ 0 & \cdots & \frac{1}{\sqrt{\lambda_N}} \end{pmatrix} \times Emat^T \tag{8}$$

$$Emat = [e_1, e_2 \ldots e_N]$$

where $\lambda_i$ is the eigenvalue with the corresponding eigenvector $e_i$ as obtained from PCA. Here the axis length of the ellipse written in this form are proportional to the square root of the diagonal elements of the diagonal matrix in the middle. It is also important to mention here that there is no analytical relationship between the best axis lengths for this ellipsoid and the $\lambda_i$s which merely give the amount of variance captured in each direction. Fig. 3 compares the $p(z)$ contours obtained using $\sqrt{\lambda_i}$ with those obtained using $\lambda_i$. The one using $\lambda_i$s is observed to have much higher eccentricity. Therefore, if such a $p(z)$ is expanded, it will increase the chances of the disturbance trajectory exiting the Lyapunov estimate in a direction transversal to the major axis, eventually resulting in a conservative CCT estimate.

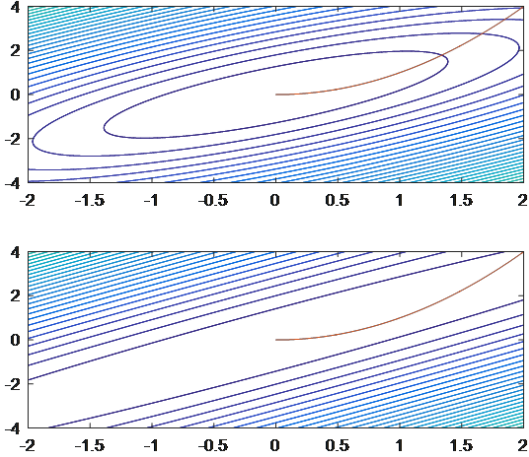

Fig. 3: $p(z)$ contours using $\sqrt{\lambda_i}$ (upper) vs. $\lambda_i$ (lower)

## IV. RESULTS

We will test our methodology on a standard 3 machine system for different scenarios and demonstrate the improvement in CCT estimates using the proposed initial choice of $p(z)$. Machine 3 is taken as reference and the system is represented in single machine reference frame. Uniform damping value $\frac{D}{M}$ of 2 is used.

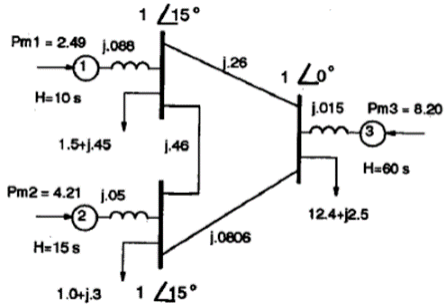

Fig. 4: 3-Machine Test Case

### A. Considerable Improvement Case

In this case, a fault is applied on bus 2 end of line 1-2. The fault is based on prior knowledge of the Lyapunov estimate shape which had a shorter axis in the direction of machine 2 states $(\delta_2, \omega_2)$. The SEP of the post-fault system is (0.3487,0.2070,0,0). Ideally, the portion of the fault trajectory data lying inside the actual SR should be used to estimate the starting value of $p(z)$. But since it's not known beforehand, we used the data from the fault sustained for 1 second. The fault trajectory vs. time is plotted in terms of the original system states $x$ and polynomial states $z$ in Fig. 5.

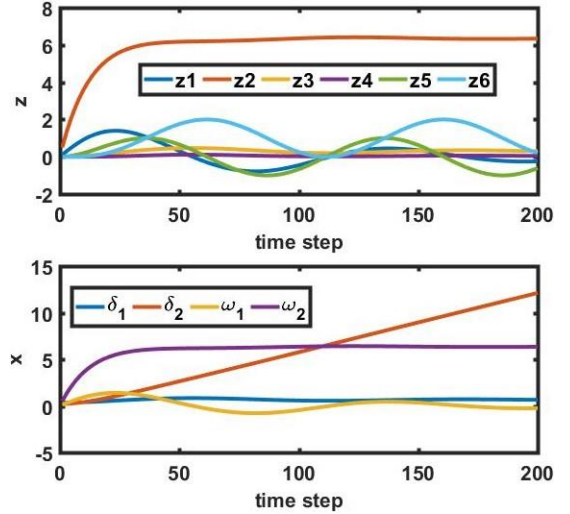

Fig. 5: Case 1 Fault Trajectory

PCA on $z$ yields the following results:
$$\lambda = [37.83, 0.7965, 0.4246, 0.0769, 0.0024, 2.2e-5]$$

$$Emat = \begin{bmatrix} 0.0183 & 0.6410 & -0.0193 & 0.7669 & 0.0143 & 0.0016 \\ 0.9840 & 0.0148 & -0.1703 & -0.0396 & -0.0312 & 0.0075 \\ 0.0494 & -0.0219 & 0.0961 & 0.0026 & 0.9409 & -0.3203 \\ 0.0086 & -0.0066 & 0.0330 & -0.0018 & 0.3185 & 0.9473 \\ 0.0026 & 0.7503 & 0.2258 & -0.6214 & -0.0047 & -0.0022 \\ 0.1702 & -0.1596 & 0.9536 & 0.1554 & -0.1100 & 0.0014 \end{bmatrix}$$

The $A$ matrix can be found using the discussion in Section III.B. The final SR estimates obtained using an identity matrix (spherical $p(z)$) for $A$ and the proposed matrix are shown in Fig. 6 below. From the figure it becomes clear that the two shapes differ considerably with the proposed approach covering more area along the $\delta_2$ axis.

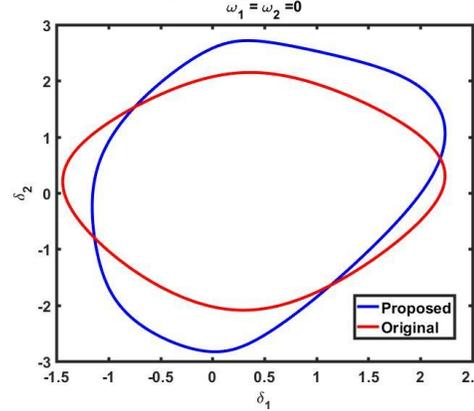

Fig. 6: Case 1 SR estimates

Let us now plot the value of the two Lyapunov functions along the fault trajectory. The CCT estimate for this fault using the proposed methodology is 0.44s while the original estimate was 0.37s. The larger CCT indicates a significant reduction in the conservativeness of the estimate.

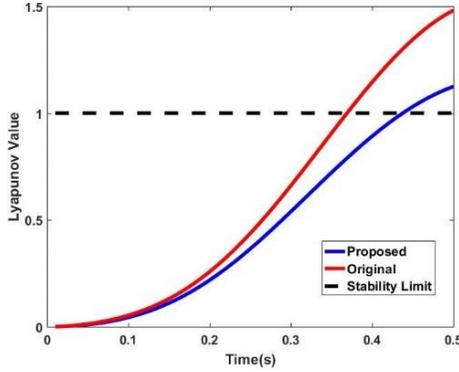

**Fig. 7: Case 1 Lyapunov function value along fault trajectory**

*B. Marginal Benefit Case*

In this case, a fault is applied on bus 1 end of line 1-2. As the SR estimate with an identity *A* matrix was significantly bigger in the $\delta_1$ direction, the proposed technique is expected to give similar results as the original with respect to the CCT estimate. The SR estimates that were obtained are shown in Fig. 8, which confirm that using the proposed technique there is only a marginal expansion of the SR estimate along the $\delta_1$ axis. Plotting the Lyapunov values along the fault trajectory we can see that there is no improvement in the CCT estimate (Fig. 9). This can be attributed to the fact that the expanding interior algorithm was able to make the original Lyapunov estimate intersect the stability boundary in the same direction with no further scope for improvement.

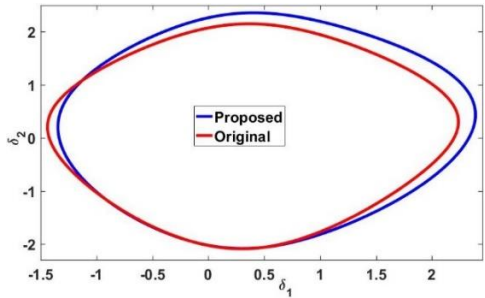

**Fig. 8: Case 2 SR estimates**

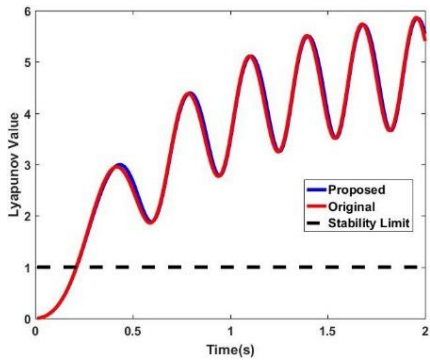

**Fig. 9: Case 2 Lyapunov function value along fault trajectory**

## V. CONCLUSIONS AND FUTURE WORK

In this paper, a methodology for estimating the relevant portion of SR for classical power system models is proposed using Lyapunov's second method and SOS programming. We first start by demonstrating the impact of the shape of the starting expanding region $P_\beta$ on the final SR estimate. Then, a PCA based approach is used that uses the disturbance trajectory data to derive a suitable ellipsoidal shape. It was seen from the results obtained for a 3 machine test system that when compared to a spherical $P_\beta$, the proposed approach considerably reduced the conservativeness in the CCT estimate. While the results are satisfactory, there are a few challenges that must be addressed. The biggest challenge is the computational limitations when dealing with large scale power systems. In this regard, the idea of decomposition using a vector Lyapunov function (found to be applicable to large interconnected systems in [11]) can be explored. Needless to say, different $P_\beta$ shapes should be used for each subsystem's SR estimate. Also, while the PCA approach is good enough for quadratic Lyapunov functions, the question of modifying the approach for higher degrees is yet to be answered.

## VI. REFERENCES

[1] H. D. Chiang, "Direct methods for stability analysis of electric power systems: theoretical foundation," *BCU Methodologies and Applications,* 1st Ed., Wiley, 2011.

[2] A. Llamas, J. D. L. R. Lopez, L. Mili, A. G. Phadke, and J. S. Thorp, "Clarifications of the BCU method for transient stability analysis," *IEEE Trans. Power Syst.*, vol. 10, no. 1, pp. 210-219, Feb. 1995.

[3] A. M. Lyapunov, *General Problem of the Stability of Motion*. London, Washington, DC: CRC Press, 1992.

[4] V. I. Zubov, *Methods of A.M. Lyapunov and their application*. Izdatel'stvo Leningradskogo Universiteta, 1961.

[5] Y. N Yu, and K. Vongsuriya, "Nonlinear power system stability study by Liapunov function and Zubov's method," *IEEE Trans. Power App. Syst.*, vol. PAS-86, no. 12, pp. 1480–1485, Dec. 1967.

[6] P. A. Parrilo, "Structured semidefinite programs and semialgebraic geometry methods in robustness and optimization," Ph.D. Dissertation, California Institute of Technology, CA, 2000.

[7] A. Papachristodoulou, and S. Prajna, "A tutorial on sum of squares techniques for systems analysis," in *Proc. American Control Conf.*, Portland, OR, vol. 4, pp. 2686-2700, June. 2005.

[8] Z. W. Jarvis-Wloszek, "Lyapunov based analysis and controller synthesis for polynomial systems using sum-of-squares optimization," Ph.D. Dissertation, University of California, Berkley, 2013.

[9] M. Anghel, F. Milano, and A. Papachristodoulou, "Algorithmic construction of lyapunov functions for power system stability analysis," *IEEE Trans. Circuits Syst. I, Reg. Papers*, vol. 60, no. 9, pp. 2533-2546, Sep. 2013.

[10] S. Kundu and M. Anghel, "Stability and Control of Power Systems using Vector Lyapunov Functions and Sum-of-Squares Methods," *European Control Conference (ECC)*, Linz, pp. 253-259, Mar. 2015.

[11] J. Anderson and A. Papachristodoulou, "A decomposition technique for nonlinear dynamical system analysis," *IEEE Trans. Autom. Control*, vol. 57, no. 6, pp. 1516–1521, Jun. 2012.

[12] C. Mishra, J. S. Thorp, V. A. Centeno, and A. Pal, "Stability region estimation under low voltage ride through constraints using sum of squares," accepted for publication in *IEEE North American power Symp. (NAPS)*, Morgantown, WV, 2017.

[13] S. Wold, K. Esbensen, and P. Geladi, "Principal Component Analysis," *Chemom. Intell. Lab. Syst.*, vol. 2, pp. 37-52, Aug. 1987.